# STLCutters.jl: A Scalable Geometrical Framework Library for Unfitted Finite Element Discretisations

Pere A. Martorell[1,*] and Santiago Badia[2]

Abstract. Approximating partial differential equations for extensive industrial and scientific applications requires leveraging the power of modern high-performance computing. In large-scale parallel computations, the geometrical discretisation rapidly becomes a bottleneck in the simulation pipeline. Unstructured mesh generation is hardly automatic, and meshing algorithms cannot efficiently exploit distributed-memory computers. Besides, partitioning of unstructured meshes relies on graph partitioning strategies, which scale poorly. As a result, the use of dynamic load balancing for locally refined meshes becomes prohibitive. Adaptive Cartesian meshes are far more advantageous, providing cheap and scalable mesh generation, partitioning, and balancing compared to unstructured meshes. However, Cartesian meshes are not suitable for complex geometries when using standard discretisation techniques. Unfitted finite element methods are a promising solution to the abovementioned problems. These numerical schemes rely on Cartesian meshes and can handle complex geometries. Nevertheless, their application is usually constrained to implicit (level set) geometrical representations. The extension to general geometries, e.g., provided by an STL surface mesh, requires advanced intersection algorithms. This work presents an efficient parallel implementation of all the geometric tools required, e.g., for unfitted finite element methods (in a broad sense), for explicit boundary representations. Such geometries can readily be generated using standard computer-aided design tools. The proposed geometrical workflow utilise a multilevel approach to overlapping computations, effectively eliminating bottlenecks in large-scale computations. The numerical results demonstrate perfect weak scalability over 13,000 processors and one billion cells. All these algorithms are implemented in the open-sopurce STLCutters.jl library, written in the Julia programming language. The library is designed to be used in conjunction with the Gridap.jl library provides a high-level interface to the finite element method.

## 1. Introduction

Nowadays, the most complex scientific simulations can only be achievable through parallel algorithms on distributed-memory machines. Large-scale parallel computations involve efficient communications to exploit the power of current high-performance computing (HPC) resources. Many finite element (FE) methods efficiently approximate the solution of partial differential equations (PDEs) in parallel. The balancing domain decomposition by constraints (BDDC) [1] and algebraic multi-grid (AMG) [2, 3] methods are popular solvers that have demonstrated high parallel scalability.

The discretisation of PDEs using FE at large scales involves distributed-memory simulations. For non-trivial geometries, the geometrical discretisation relies on generating unstructured body-fitted meshes, a challenging task usually involving human intervention and preventing an automatic simulation workflow [4]. Enforcing mesh conformity prevents the efficient exploitation of supercomputers, and mesh partitioning algorithms rely on graph partitioning techniques [5, 6] that are not scalable. Thus, the geometrical discretisation rapidly becomes a bottleneck, even a showstopper, of the simulation pipeline for large-scale simulations of problems on complex geometries. The situation is even more dramatic when adaptive mesh refinement (AMR) is required (e.g., in multiscale and multiphysics simulations) since dynamic load-balancing is unaffordable due to the performance overhead of graph-partitioning strategies.

The situation is much simpler when problems are posed in simple geometries (e.g., cubes). Locally adapted Cartesian meshes can exploit scalable and dimension-agnostic mesh generators and partitioners [7, 8]. The generation and load-balancing of octree-meshes is efficient thanks to space-filling curve techniques [9]; see, e.g., the highly scalable p4est framework [10] for handling forests of octrees on hundreds of thousands of processors.

[1] CIMNE, Centre Internacional de Mètodes Numèrics a l'Enginyeria, Campus Nord, 08034 Barcelona, Spain.
[2] School of Mathematics, Monash University, Clayton, Victoria, 3800, Australia.
* Corresponding author.
E-mails: pere.antoni.martorell@upc.edu (PM) santiago.badia@monash.edu (SB).

Unfitted (also known as embedded or immersed) FE methods [11] can overcome some of the limitations of body-fitted meshes. Unfitted methods combine a (e.g., Cartesian) background mesh for functional discretisation with quadrature rules to integrate the intersection between cells in the background mesh and the domain interior. However, naive unfitted discretisations may lead to unstable and severe ill-conditioned discrete problems [11] unless a specific technique mitigates the problem. The size and aspect ratio of the intersection of a background cell and the physical domain are not bounded. Despite the vast literature on the topic, unfitted FE formulations that solve these issues are quite recent. Stabilised formulations based on the so-called *ghost penalty* were originally proposed in [12] for Lagrangian continuous FEs and has been widely used since [13]. The so-called *cell aggregation* or *cell agglomeration* techniques are an alternative way to ensure robustness concerning cut location. This approach is very natural in discontinuous Galerkin (DG) methods, as their formulation on agglomerated meshes is straighforward [14]. These techniques present extensions with $C^0$ Lagrangian FE in [15] and mixed methods in [16]. These unfitted formulations enjoy good numerical properties, such as condition number bounds, stability, optimal convergence, and continuity with respect to data. Distributed-memory implementations for large-scale problems have been designed [17], and error-driven $h$-adaptivity and parallel tree-based meshes have also been exploited [8]. Another challenge of unfitted FE methods is the generation of meshes for wall-bounded turbulent flows, which require anistropic meshes aligned with the wall. Such situations could be handled by combining unfitted methods with overlapping mesh strategies (a.k.a. Chimera methods) [18].

The generation of *body-fitted quadratures* required in unfitted formulations is more amenable to parallelisation than the generation of the body-fitted meshes required by standard methods. Most ingredients in the computation of these body-fitted quadratures are cell-local and thus embarrassingly parallel. However, most proposed methods rely on implicit geometry representations, e.g., level sets, reducing their industrial application [8]. Recently, in [19], the authors have presented robust schemes for explicit boundary representation (BREP), which have been extended in [20] to high order BREPs, e.g., computer-aided design (CAD) models. Unlike level set methods, explicit geometry representations require global operations to define the inside and outside domains, e.g., ray-tracing and propagation techniques. These operations are not trivial to parallelise and can become a bottleneck in large-scale simulations if not computed efficiently.

This work is the last step of a research project started several years ago, with the aim to provide a scalable and automatic geometrical framework for unfitted discretisations of PDEs. In [19] and [20], the authors presented robust algorithms for the computation of body-fitted quadratures for unfitted formulations for explicit BREP geometries and high-order CAD models, respectively. In this work, we present a scalable geometrical workflow for explicit (and implicit) geometry representations. It combines the cell-local algorithms in [19], which are embarrassingly parallel, with a novel scalable implementation of the in-out definition. Namely, we define a multilevel propagation algorithm that overlaps coarse and fine computations to achieve perfect scalability. The presented algorithm, inspired by the multilevel BDDC [21, 22], shows excellent weak scalability over 10k processors. Besides, we present readily available open-source software for all proposed algorithms.

The outline of this article is as follows. In Sec. 2, we introduce the parallel unfitted FE methods used in this work. We present the distributed parallel geometrical framework for unfitted discretisations in Sec. 3. In Sec. 4, we present STLCutters.jl, the open-source software that provides a scalable implementation of the proposed geometrical workflow. We analyse the framework's scalability on distributed-memory platforms in Sec. 5. Finally, in Sec. 6, we summarise the main findings of this work.

## 2. Distributed unfitted finite element method

In this section, we provide a brief introduction to unfitted FE methods, in order to identify the geometrical requirements of these formulations. We focus on FE techniques in a broad sense, which includes DG or hybrid formulations. The discussion is PDE-agnostic, and can readily be applied to interface and bulk-coupling multiphysics.

Let us consider a Lipschitz domain $\Omega \subset \mathbb{R}^d$, with $d \in \{2, 3\}$ the number of spatial dimensions. We describe the domain's boundary $\partial \Omega$ with a parametric oriented surface mesh $\mathcal{B}$. We aim to solve a system of PDEs that can involve Dirichlet boundary conditions on $\Gamma_D$ and Neumann boundary conditions on $\Gamma_N$.

$\Gamma_N$ and $\Gamma_D$ are a partition of $\partial\Omega$. The geometrical representation, e.g., CAD model, must respect this partition. Therefore, we consider $\mathcal{B}_D$ and $\mathcal{B}_N$ a partition of $\mathcal{B}$, that represent $\Gamma_D$ and $\Gamma_N$ respectively.

This work aims to define an efficient method for the parallel implementation of unfitted FE discretisations generated from $\mathcal{B}$. Unfitted discretisations utilise a background mesh $\mathcal{T}^{\text{bg}}$ mesh instead of body-fitted meshes. This background mesh is an arbitrary partition of an artificial domain $\Omega^{\text{art}}$ such that $\Omega^{\text{art}} \supset \Omega$ (see Fig. 1). $\Omega^{\text{art}}$ can be as simple as a bounding box. The $\mathcal{T}^{\text{bg}}$ is a more straightforward partition than a body-fitted mesh, e.g., a Cartesian grid or a refinement of a hexahedral mesh.

In distributed-memory computation, we subdivide the domain $\Omega^{\text{art}}$ into $S$ subdomains $\Omega^{\text{art}}_s$, $s = 1, \ldots, S$ (see Fig. 1). Subdomain partition of uniform and adaptive Cartesian meshes is very efficient since one can efficiently exploit octrees using space-filling curves [9], avoiding graph partition algorithms [6]. In unfitted formulations, one can readily use weighting factors (see [10]) to balance active cells.

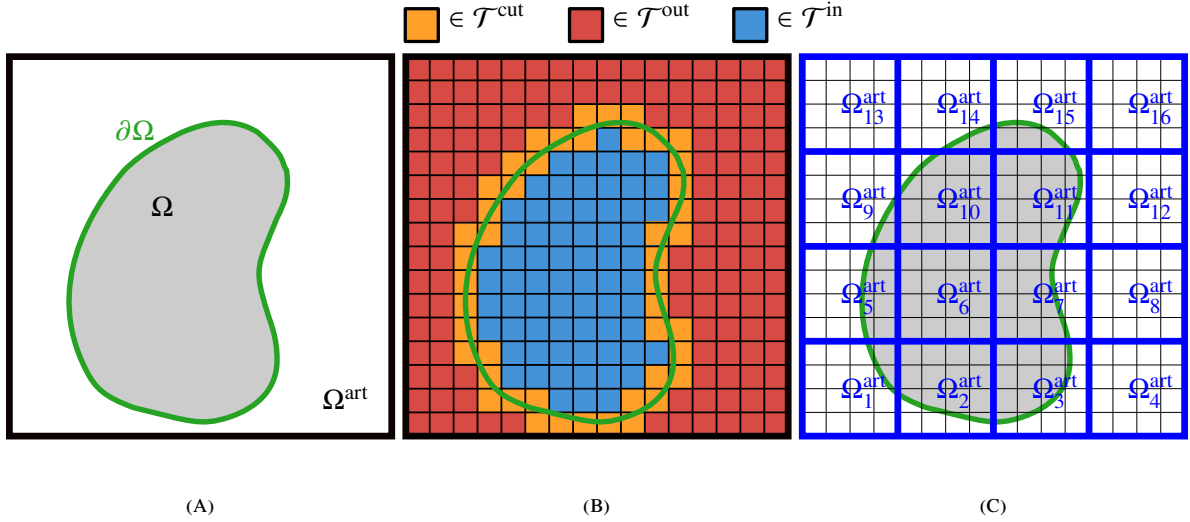

FIGURE 1. Unfitted FE representation in a distributed-memory computation. In (a), a surface mesh $\mathcal{B}$ represents the domain's boundary $\partial\Omega$. The domain $\Omega$ is embedded in $\Omega^{\text{art}}$. The artificial domain $\Omega^{\text{art}}$ is discretised in a Cartesian background mesh $\mathcal{T}^{\text{bg}}$. In (b), We classify the cells in $\mathcal{T}^{\text{bg}}$ interior, exterior, and cut cells. For parallel computations, in (c), we divide $\Omega^{\text{art}}$ in $S$ subdomains $\Omega^{\text{art}}_s$, $s = 1, \ldots, S$. The background mesh is partitioned accordingly. The partition into subdomains coincides with the coarse mesh $\mathcal{T}^{\text{coarse}}$. We consider a naive partition of the mesh for simplicity in the illustration. In practice, one should balance active cells instead of cells.

The presented unfitted FE formulation accommodates different methods from the literature, e.g., the extended FE method (XFEM) [23], the cutFEM [13], the aggregated finite element method (AgFEM) [15], or the finite cell method [24]. Furthermore, the presented algorithm is agnostic to the solver utilised for parallel unfitted FE methods, e.g., BDDC [7] or AMG [2].

To solve PDEs on the unfitted discretisation, we need to classify the background cells into interior, exterior, and cut, $\mathcal{T}^*$, $* \in \{\text{in}, \text{out}, \text{cut}\}$, respectively (see Fig. 1). The functional discretisation does not consider the exterior cells $\mathcal{T}^{\text{out}}$. Thus, we consider the active mesh $\mathcal{T} \doteq \mathcal{T}^{\text{bg}} \setminus \mathcal{T}^{\text{out}}$ for the FE discretisation. The unfitted FE techniques utilise standard FE spaces on $\mathcal{T}$, $V$, to solve and test the weak form of the PDEs. The unfitted problem reads as follows: find $u \in V$ such that

$$a(u, v) = l(v), \quad \forall v \in V,$$

where

$$a(u, v) = \int_\Omega L_\Omega(u, v) d\Omega + \int_{\Gamma_D} L_D(u, v) d\Omega + \int_\mathcal{F} L_{sk}(u, v) d\Gamma,$$

and

$$l(v) = \int_\Omega F_\Omega(v) d\Omega + \int_{\Gamma_N} F_N(v) d\Gamma + \int_{\Gamma_D} F_D(v) d\Gamma.$$

Here, we include the differential operator, the source term, and additional stabilisation terms within the bulk terms $F_\Omega$ and $L_\Omega$. The integration of the operations $F_D$ and $L_D$ on $\Gamma_D$ are related to the weak imposition of the Dirichlet boundary conditions, commonly utilising Nische's method. The term $F_N$ integrated on $\Gamma_N$ impose the Neumann boundary conditions. The skeleton $\mathcal{F}$ represents the interior faces of the active mesh $\mathcal{T}$. Integrating $L_{sk}$ on $\mathcal{F}$, we introduce additional penalty terms, e.g., ghost penalty stabilisation techniques or weak imposition of continuity in DG methods.

In FE methods, the integration of *piecewise* polynomials relies on a cell-wise decomposition of bulk and surface terms. In unfitted FE methods, the integration is defined only in the domain's interior. Bulk integration meshes $\mathcal{T}^{\text{int}} \doteq \{K \cap \Omega : K \in \mathcal{T}\} = \mathcal{T}^{\text{in}} \cup \mathcal{T}^{\text{clip}}$, where the clipped mesh reads as $\mathcal{T}^{\text{clip}} \doteq \{K \cap \Omega : K \in \mathcal{T}^{\text{cut}}\}$. A cell-wise geometrical algorithm, e.g., algorithms in [19] and [20], performs the intersection of the cut cells $K \in \mathcal{T}^{\text{cut}}$. However, the definition of the interior cells $K \in \mathcal{T}^{\text{in}}$ is not trivial in distributed memory computations. The propagation of inside cells can become a bottleneck in large-scale computations. The integration of the surface terms requires the intersection of the boundary with the background mesh $\mathcal{B}^{\text{int}} \doteq \{\mathcal{B} \cap K : K \in \mathcal{T}^{\text{cut}}\}$.

The parallel unfitted FE methods each background cell $K \in \mathcal{T}^{\text{bg}}$ belongs to a single task $s \in 1, ..., S$ (see Fig. 1c). Therefore, one can define the above meshes in each of the $S$ tasks, e.g., $\mathcal{T}_s^*$, $* \in \{\text{in, out, cut, bg, int}\}$, $\mathcal{B}_s$ and $\mathcal{B}_s^{\text{int}}$. In the proposed algorithms, each coarse cell $K_s^{\text{coarse}} \in \mathcal{T}^{\text{coarse}}$ represents an artificial subdomain $\Omega_s^{\text{art}}$. Unlike distributed FE solvers, we do not require *ghost* cells (i.e., cells that belong to other processors) for the geometrical discretisation.

## 3. Distributed intersection algorithm

In this section, we define an algorithm that intersects the background distributed cells $K \in \mathcal{T}^{\text{bg}}$ and a physical domain $\Omega$ described by an oriented surface mesh $\mathcal{B}$. First, in Sec. 3.1, we briefly expose the algorithms for the bulk intersections $K \cap \Omega$ and the boundary intersections $K \cap \mathcal{B}$ proposed in [19]. Then, in Sec. 3.2, we propose a propagation method for classifying non-intersected cells $K \cap \partial\Omega = \emptyset$ as interior or exterior cells. Finally, in Sec. 3.3, we design a global algorithm that overlaps classifications and intersections in a two-level distributed mesh.

3.1. **Local intersection.** The unfitted FE methods' main complication resides in the cut cells' integration. For this purpose, we need to compute the intersections of the background cells $K \in \mathcal{T}^{\text{bg}}$ with the domain $\Omega$. Each integration method relies on a different underlying data structure to store the resulting intersection $\mathcal{T}^{\text{clip}}$. In implementing STLCutters [25], we follow the intersection methods described in [19], which generates a two-level simplex mesh.

In [19], we described an intersection algorithm for embedded discretisations on piecewise linear explicit boundary representations. Such algorithms have proven robust and accurate up to machine precision. The algorithm is based on the intersection of convex polytopes with the half-spaces of the boundary planes; see [26, 27]. The approach in [19] performs a convex decomposition step to deal with general non-convex geometries. The convex decomposition consists on recursively splitting a non-convex and a background cell using the planes passing through *reflex edges* (the edges that connect two non-convex faces). The result is a set of convex surface components and a set of background pieces. We can clip the background cell pieces with the planes of the corresponding convex surface component. It generates a set of convex polytopes that represent the intersection of the background cell with interior of the domain.

However, one can utilise meshes of general polytopes or combine moment fitting techniques with the Stokes theorem [28, 29]. For non-linear parametrisations of $\mathcal{B}$, we can consider the high-order intersection algorithms in [20].

3.2. **Local classification.** The definition of the non-intersected background cells $K \in \mathcal{T}^{\text{bg}} \setminus \mathcal{T}^{\text{cut}}$ is a simple task in serial computations. Taking advantage of the information the intersected cells give, $K \in \mathcal{T}^{\text{cut}}$, a depth-first traversal propagation is an efficient choice. In contrast to ray tracing algorithms [30], this method does not require floating point operations.

The algorithm utilises nodal propagation of the relative positions, i.e., interior and exterior. Nodal propagation is trivial in conformal meshes, e.g., Cartesian meshes. However, we may need to propagate across the hanging nodes in non-conformal meshes, e.g., octree meshes from p4est. We do not need to consider the hanging nodes in 2:1 $k$-balanced octrees, i.e., a maximum of one hanging node per edge (see [10]).

3.3. **Global distributed algorithm.** Each task in a distributed-memory environment computes local computations in a subdomain $\Omega_s^{\text{art}} \subset \Omega^{\text{art}}$. This subdomain is discretised with a local background mesh $\mathcal{T}_s^{\text{bg}}$ in each task. In large-scale computations, a single task can not contain the global mesh $\mathcal{T}^{\text{bg}}$. One additional task contains a coarse background mesh $\mathcal{T}^{\text{coarse}}$. Each cell of this coarse mesh $K_s^{\text{coarse}} \in \mathcal{T}^{\text{coarse}}$ represents a subdomain $\Omega_s^{\text{art}}$, thus $K_s^{\text{coarse}}$ and $\Omega_s^{\text{art}}$ belong to a single task $s$. This coarse mesh can be reused in the solver stage, e.g., in the BDDC solver [21]. Similarly to BDDC, one can consider a multilevel coarse mesh to handle a larger number of tasks. However, in this work, we consider a single coarse mesh (see Fig. 2). Alternatively, one could consider scalable sparse communications protocols [31].

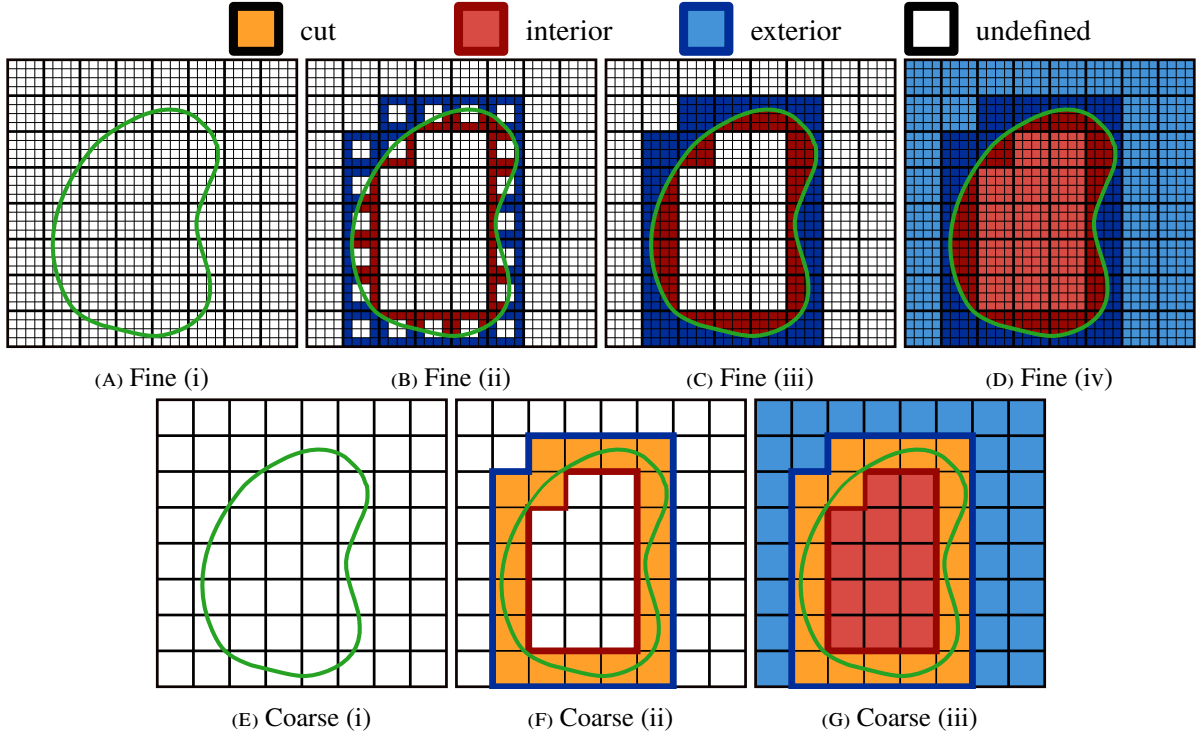

FIGURE 2. Representation of the two-level propagation algorithm. The domain boundary $\partial \Omega$ intersects the fine (a) and coarse mesh (e). First, we intersect the domain boundaries $\Omega_s^{\text{art}}$ (b) to define the coarse cells $K_s^{\text{coarse}} \in \mathcal{T}^{\text{coarse}}$ (f). The `gather` message-passing interface (MPI) command performs the fine to coarse communications. Then, we propagate the coarse cells (g) while intersecting the rest of the fine subdomains (c). Finally, after a coarse to fine communication (`scatter`), we define the cells in the untouched subdomains (d).

The main advantage of using a two-level distributed algorithm is the reduction of communications. The centralised communications go through the root processor that contains the coarse mesh. The processors communicate with the root processor using `gather` and `scatter` MPI commands. The processors can also communicate with the nearest neighbours using `send`, `recv`, or `sendrecv` MPI commands. The nearest-neighbour communications enforce consistency in the elements near the interface between subdomains. We perform the fewest nearest-neighbor communications to maintain the algorithm's scalability.

Alg. 1 performs the intersection of the background cells and the inside propagation in a distributed-memory computation. We developed an efficient and scalable algorithm by overlapping computations. Each MPI task runs the algorithm. The fine tasks execute the code blocks within $s \neq s^{\text{coarse}}$ (see line 1, 12 and 20) and the coarse task runs the blocks in $s = s^{\text{coarse}}$ (see line 15). All the tasks execute the code outside these conditions.

The first part of the algorithm, line 2-3, performs the intersection of the cells near the boundary of each subdomain $\partial \Omega_s^{\text{art}}$ (or $\partial K_s^{\text{coarse}}$); see Fig. 2b. Here, we define the location map of the faces as $\mathcal{L}: F \mapsto \{\text{in}, \text{out}, \text{cut}\}$ for $F \in \text{faces}(\mathcal{T})$. In line 4, we extract the location map $\mathcal{L}_s$ from the clipped boundary mesh $\mathcal{T}_s^{\text{clip,bnd}}$, at this point only the faces touching the clipped cells have a defined location.

**Algorithm 1** `distributed_intersection`($\mathcal{T}_s^{\text{bg}}, K_s^{\text{coarse}}, \Omega_s$)

1: **if** $s \neq s^{\text{coarse}}$ **then**
2:    $\mathcal{T}_s^{\text{bnd}} \leftarrow \{K \in \mathcal{T}_s : K \cap \partial K_s^{\text{coarse}} \neq \emptyset\}$
3:    $\mathcal{T}_s^{\text{clip,bnd}} \leftarrow \{K \cap \Omega_s : K \in \mathcal{T}_s^{\text{bnd}}, K \cap \partial \Omega_s \neq \emptyset\}$
4:    $\mathcal{L}_s \leftarrow \texttt{location\_map}(\mathcal{T}_s^{\text{clip,bnd}})$
5:    $\mathcal{L}_s \leftarrow \texttt{propagate\_location}(\mathcal{T}_s^{\text{bnd}}, \mathcal{L}_s)$
6:    $\mathcal{L}_s \leftarrow \texttt{sendrecv}(\mathcal{L}_s)$
7:    $\mathcal{L}_s^{\text{coarse}} \leftarrow \{(F, \mathcal{L}_s(\mathcal{F}(F))) : F \in K_s^{\text{coarse}}\}$
8: **end if**
9: $\mathcal{L}^{\text{coarse}} \leftarrow \texttt{gather}(\mathcal{L}_s^{\text{coarse}})$
10: **if** $s \neq s^{\text{coarse}}$ **then**
11:    $\mathcal{T}_s^{\text{bulk}} \leftarrow \mathcal{T}_s \setminus \mathcal{T}_s^{\text{bnd}}$
12:    $\mathcal{T}_s^{\text{clip}} \leftarrow \mathcal{T}_s^{\text{clip,bnd}} \cup \{K \cap \Omega_s : K \in \mathcal{T}_s^{\text{bulk}}, K \cap \partial \Omega_s \neq \emptyset\}$
13:    $\mathcal{L}_s \leftarrow \texttt{location\_map}(\mathcal{T}_s^{\text{cut}})$
14:    $\mathcal{L}_s \leftarrow \texttt{propagate\_location}(\mathcal{T}_s^{\text{bg}}, \mathcal{L}_s)$
15: **else**
16:    $\mathcal{T}^{\text{coarse}} \leftarrow \mathcal{T}_s$
17:    $\mathcal{L}^{\text{coarse}} \leftarrow \texttt{propagate\_location}(\mathcal{T}^{\text{coarse}}, \mathcal{L}^{\text{coarse}})$
18: **end if**
19: $\mathcal{L}_s^{\text{coarse}} \leftarrow \texttt{scatter}(\mathcal{L}^{\text{coarse}})$
20: **if** $s \neq s^{\text{coarse}}$ **then**
21:    **if** $\mathcal{T}_s^{\text{clip}} = \emptyset$ **then**
22:      $\mathcal{L}_s \leftarrow \{(F, \mathcal{L}_s^{\text{coarse}}(K_s^{\text{coarse}})) : F \in \mathcal{T}_s\}$
23:    **end if**
24:    $\mathcal{T}_s^{\text{in}} \leftarrow \{K \in \mathcal{T}_s^{\text{bg}} : \mathcal{L}_s(K) = \text{in}\}$
25: **end if**
26: **return** $\mathcal{T}_s^{\text{in}} \cup \mathcal{T}_s^{\text{clip}}$

The map $\mathcal{L}_s$ is indexed with the local faces of $\mathcal{T}_s^{\text{bg}}$. Then, in line 5, we propagate the location $\mathcal{L}_s$ through the domain boundary cells $\mathcal{T}_s^{\text{bnd}}$. The propagation starts from the defined cells $\mathcal{T}_s^{\text{clip,bnd}}$ as described in Sec. 3.1. The line 6 enforces consistency of the interface faces with nearest-neighbor communications, e.g., `sendrecv` MPI command. At this point, the location of the subdomains is defined with minimum operations and communications. Thus, we build a map for the location of the local coarse entities $\mathcal{L}_s^{\text{coarse}}$ in line 7. For non-intersected subdomains, i.e., $\mathcal{T}_s^{\text{clip}} = \emptyset$, the location of the entire subdomain is set as *undefined*. The coarse task collects the local locations of the coarse entities in line 9 with a `gather` communication (see Fig. 2f).

In the second part, we overlap the intersection of the rest of the fine cells (line 12) and the classification of the coarse mesh (line 15). This overlapping is crucial to avoid the idling of the main bunch of processors. The intersection of the bulk cells $\mathcal{T}_s^{\text{bulk}}$ in line 12 is analogous to the intersections of the boundary cells $\mathcal{T}_s^{\text{bnd}}$ in line 2. Here, we intersect and classify the entire subdomains. Classifying the coarse cells $\mathcal{T}^{\text{coarse}}$ in line 17 utilises serial propagation algorithms. After this part, the coarse sends the coarse location to each processor through the `scatter` command (line 19). Finally, we define the local cells in the untouched domains (line 22). This algorithm returns a mesh ready for integration $\mathcal{T}^{\text{in}} \cup \mathcal{T}^{\text{clip}}$.

It is important to note that in Fig. 2, some subdomains are entirely outside and have empty contributions to the system. These subdomains do not present an issue for interface problems that have contributions on both interior and exterior parts. However, for boundary value problems, one may need to redistribute the mesh to balance the load. We present the redistribution techniques in the next section.

## 4. STLCutters.jl

This section introduces the open-source software STLCutters.jl, which includes an efficient parallel implementation of the geometrical workflow defined in the previous section. STLCutters.jl [25] is fully written in Julia [32]. In the design of the package, we take advantage of the features of the language, e.g., multiple dispatch, multi paradigm, just-in-time (JIT) compilation and high-level and low-level programming. Additionally, it utilises the rich and growing Julia ecosystem of packages through the built-in package manager. In particular, STLCutters strongly depends on Gridap [33, 34], a FE library with a high-level user interface that resembles mathematical notation, and the Gridap ecosystem, e.g., GridapEmbedded [35] (for unfitted discretisations) and GridapDistributed [36] (for the parallel implementation). As far as we know, most FE packages in the literature do not provide an implementation of unfitted FE methods, since these schemes have specific requirements (e.g., geometrical tools and runtime definition of quadratures on cut cells) not common in standard finite element method (FEM). A notable exception is ngsxfem, an add-on to the NGSolve package [37].

STLCutters utilises data structures and procedures from Gridap and a similar coding style. Thus, like in Gridap, STLCutters has a clear and intuitive high-level application programming interface (API). Indeed, STLCutters API is an extension of GridapEmbedded API. We demonstrate such simplicity in Lst. 1 with an example of the Poisson equation on an embedded distributed domain. We note that Gridap is a general PDE solver that is not bond to specific applications. It has been applied to complex physical systems, e.g., inductionless magnetohydrodinamics [38], the Keller-Segel equations [39], or fluid-structure interaction on moving domains [40], to mention a few. However, we consider the Poisson problem here for simplicity in the exposition, since the emphasis is in the geometrical workflow.

In the example, lines 12 and 15 represent the interface of STLCutters. We load the geometry from an STL file in line 12. In line 13, we compute a bounding box for this geometry. We provide the desired number of subdomain partitions per direction in line 6 (see Fig. 3a) and the number of cells per direction in line 7. We create a background Cartesian mesh in line 14 with all this information. In line 15, we cut the geometry with the background mesh using the algorithms presented in Sec. 3.

The rest of the code is a standard FE driver with Gridap and GridapEmbedded [41], which we briefly describe here for completeness. In this example, to fix the so-called small cut cell problem, we use the aggregated finite element method (see [15]), which requires the aggregation of cut cells to interior cells, computed in line 17. With these aggregates, we can define the aggregated Lagrangian FE spaces in lines 25-27, which we will use in the unfitted discretisation. We extract the active domain (interior part) $\Omega$ and the embedded boundary $\Gamma$ in lines 19 and 20. We define the normal vector to the boundary in line 21. The *measure* in line 22 represents the body-fitted quadrature required to integrate the bulk terms. The boundary *measure* is computed in line 23. With these quadratures, we can define the (bi)linear forms in an unfitted FE discretisation with a weak imposition of Dirichlet boundary conditions with Nitsche's method in lines 33-38. Finally, we solve the problem in line 41 and write the results in a VTK file in line 42.

The parallel implementation of the STLCutters relies on GridapDistributed and PartitionedArrays. GridapDistributed offers a distributed extension of the Gridap types and procedures. PartitionedArrays provides a layer of abstraction over the MPI primitives. It is important to note that GridapDistributed also depends on PartitionedArrays. To run a serialised distributed code, one can change line 9 of Lst. 1 by `with_debug`. This feature allows debugging without MPI executions.

Adding a few of the lines of Lst. 2 into Lst. 1, the code runs with a balanced `p4est` background mesh. The `redistribute` procedure distributes the background cells with a given weight. In unfitted FEM, the exterior cells are not contributing to the FE approximation. Thus, setting zero-weight to these cells would improve the load-balancing across subdomains. One can see the differences in load balancing in Fig. 3.

In addition to load balancing, one can set AMR on the cut cells to capture the physical phenomena around the boundary. The code in Lst. 3 adapts the background model around the boundary. Since the classification of the cells is nodal-based, we need to consider 2:1 $k$-balance in the AMR; see Sec. 3.2. We remind that STLCutters captures the exact geometrical features. Unlike in level set methods, AMR will not change the geometrical description.

## 5. Numerical experiments

5.1. **Experimental setup.** The numerical experiments have been performed on Gadi, a high-end supercomputer at the NCI (Canberra, Australia) with 4962 nodes, 2074 of them powered by a 2 x 24 core Intel Xenon Platinum 8274 (Cascade Lake) at 3.2 GHz and 192 GB RAM. The algorithms presented in this work are implemented in Julia [32]. We used Julia v1.9, Intel MPI 2021.10, and STLCutters v0.2.1 during the experiments.

**Listing 1** Code usage for solving the Poisson equation on an embedded domain (293137 STL file) using the AgFEM method. This example should be run with MPI as `mpiexec -np 4 julia example.jl`.

```julia
using STLCutters
using Gridap
using GridapEmbedded
using GridapDistributed
using PartitionedArrays
parts = (1,1,4)
cells = (8,8,8)
filename = "293137.stl"
with_mpi() do distribute
  ranks = distribute(LinearIndices((prod(parts),)))
  # Domain and discretisation
  geo = STLGeometry(filename)
  pmin,pmax = get_bounding_box(geo)
  model = CartesianDiscreteModel(ranks,parts,pmin,pmax,cells)
  cutgeo = cut(model,geo)
  # Cell aggregation
  model,cutgeo,aggregates = aggregate(AggregateAllCutCells(),cutgeo)
  # Triangulations
  Ω_act = Triangulation(cutgeo,ACTIVE)
  Ω = Triangulation(cutgeo,PHYSICAL)
  Γ = EmbeddedBoundary(cutgeo)
  nΓ = get_normal_vector(Γ)
  dΩ = Measure(Ω,2)
  dΓ = Measure(Γ,2)
  # FE spaces
  Vstd = TestFESpace(Ω_act,ReferenceFE(lagrangian,Float64,1))
  V = AgFEMSpace(model,Vstd,aggregates)
  U = TrialFESpace(V)
  # Weak form
  γ = 10.0
  h = (pmax - pmin)[1] / cells[1]
  ud(x) = x[1] - x[2]
  f = 0
  a(u,v) =
    ∫( ∇(v)⋅∇(u) )dΩ +
    ∫( (γ/h)*v*u  - v*(nΓ⋅∇(u)) - (nΓ⋅∇(v))*u )dΓ
  l(v) =
    ∫( v*f )dΩ +
    ∫( (γ/h)*v*ud - (nΓ⋅∇(v))*ud )dΓ
  # Solve
  op = AffineFEOperator(a,l,U,V)
  uh = solve(op)
  writevtk(Ω,"results",cellfields=["uh"=>uh])
end
```

**Listing 2** Compute balanced model with `p4est`. It replaces line 14 in Lst. 1 before calling the `cut` function again. In order to use this code, `GridapP4est` needs to be loaded

```julia
coarse_model = CartesianDiscreteModel(pmin,pmax,cells)
model = OctreeDistributedDiscreteModel(parts,coarse_model)
cutgeo = cut(model,geo)
weights = compute_redistribute_weights(cutgeo)
model,_ = GridapDistributed.redistribute(model,weights=weights)
```

**Listing 3** AMR with `p4est`. This code should be added between lines 2 and 3 of Lst. 2

```julia
cutgeo = cut(model,geo)
flags = adapt_cut_cells(parts,cutgeo)
model,_ = Gridap.Adaptivity.adapt(model,flags)
```

5.2. **Parallel scalability.** In this section, we focus the results on the parallel performance, specifically on the algorithm's scalability. We utilise the subset of geometries from Thingi10K [42] represented in Fig. 4. These geometries are embedded in an artificial domain $\Omega^{\text{art}}$. The artificial domain is 40% larger than the bounding box of the geometries in each Cartesian direction. This value is inherited from [19] to have a significant layer of outside cells to test the propagation. The bounding box size would have the similar effect in the scalability tests than the number of background cells. The background mesh $\mathcal{T}_{bg}$ is a Cartesian partition of the artificial domain $\Omega^{\text{art}}$. We compute the embedded discretisation of the given geometry and background mesh in each data point.

We perform strong scalability tests to estimate the local size in the weak scaling tests. These tests are performed with a fixed number of cells, 3,072,000 cell background mesh ($120 \times 160 \times 160$). This

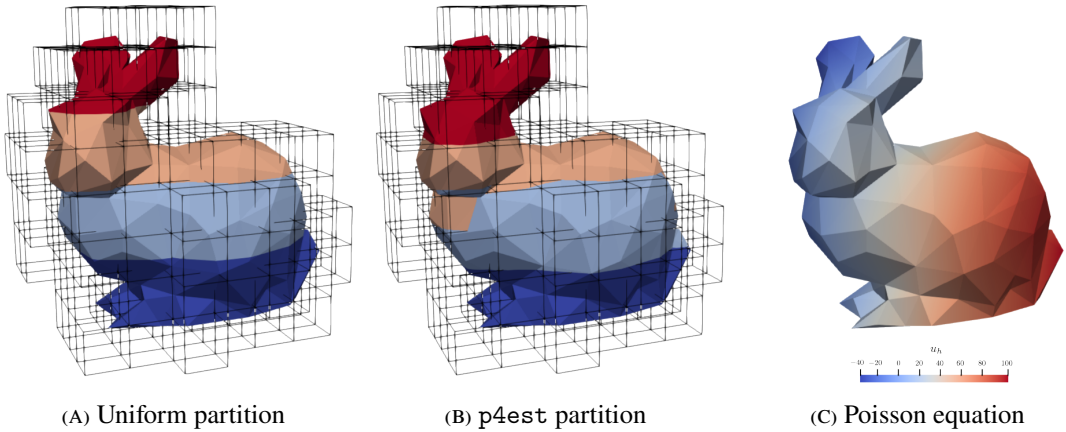

(A) Uniform partition  (B) `p4est` partition  (C) Poisson equation

FIGURE 3. Balancing differences between a uniform partition (a) $P = 1 \times 1 \times 4$ and a `p4est` partition (b). Both solve the Poisson equation (c) in a distributed domain of 4 processors in a mesh of $N = 8 \times 8 \times 8$ on `293137.stl`.

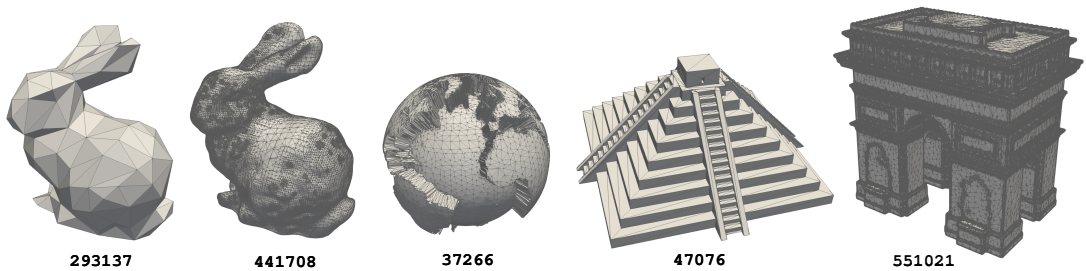

293137  441708  37266  47076  551021

FIGURE 4. Set of geometries used for testing the distributed intersection algorithms. The code of the geometries is the identifier in the Thingi10K dataset [42].

background mesh is tested with a range of central processing unit (CPU) counts, $P = \{1, 2, 4, 8, 16, 32, 48\}$. As expected, the strong scalability in Fig. 5a is not optimal. We do not aim to design a strongly scalable algorithm because cut cells are unbalanced across processors for denser Stereolithography (STL) geometries. One could readily improve the strong scalability by load-balancing the boundary cells. However, this would not be an efficient partition for FE computations, which represent the bulk of the simulation time.

The algorithm is designed to be weakly scalable. Weak scaling is essential in large-scale simulations. In the weak scaling tests, we increase the number of processors $P$ with the number of background cells $N$ while keeping the local number of cells $N/P$ constant. Based on the strong scaling tests, we test the following local sizes $N/P = \{1,000;\ 64,000\}$, i.e., $N/P = \{10^3, 40^3\}$. In Fig. 5, we show perfect weak scalability for local problems of $1k$ and $64k$ cells. Moreover, for finer STL geometries, the computational time is even decreasing with the number of processors. In these geometries, the number of faces per cut cell is still decreasing for fine background meshes. The tests are performed for $P = \{48, 384, 3072, 13824\}$ processors, i.e., $P = \{3 \times 4 \times 4;\ 6 \times 8 \times 8;\ 12 \times 16 \times 16;\ 24 \times 24 \times 24\}$. Therefore, we have tested perfect weak scalability up to 884,736,000 cells in 13,824 processors.

## 6. Conclusions

In this work, we have presented a distributed and automatic geometrical workflow for unfitted large-scale simulations on distributed-memory systems. The workflow is automatic and robust, capturing the exact geometrical features of the STL files, and has been designed to be weakly scalable. We provide the implementation of the proposed workflow in the open-source software STLCutters.jl, a Julia package that extends the Gridap.jl ecosystem. We have demonstrated the scalability of the implementation in a large-scale distributed-memory environment. The results show perfect weak scalability up to 884,736,000 cells in 13,824 processors, while the volume errors are close to machine precision for the bounding box and the interior volume. Combined with state-of-the-art FE solvers, the proposed geometrical framework, and more specifically, the STLCutters.jl implementation, is a powerful simulation pipeline for large-scale simulations in complex geometries.

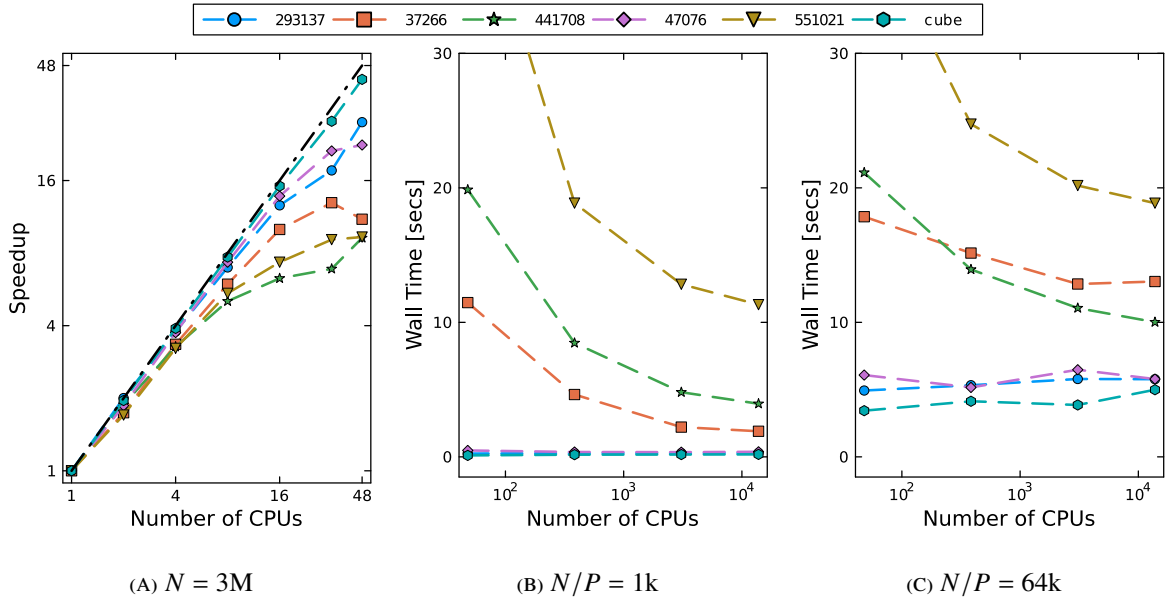

FIGURE 5. Parallel scalability of the creation of distributed embedded discretisations. The strong scaling tests in (a) are performed with a fixed background mesh of 3,072,000 cells. The weak scaling tests in (b) and (c) are computed by increasing the number of processors $P$ and the number of background cells with local loads of $N/P = \{1,000; 64,000\}$. The algorithm is not strongly scalable, but the weak scalability shows optimality up to 884,736,000 cells in 13,824 processors.

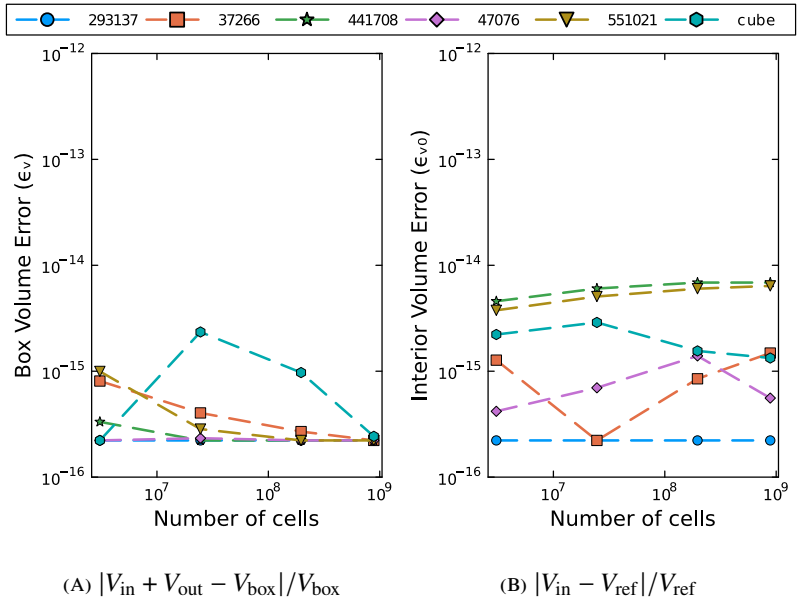

FIGURE 6. Volume errors in weak scaling test. The discretisation has bounding box volume errors (a) close to machine precision. The variation of the interior (b) volume is also low. Both cases for a local size of $N/P = 64,000$ cells.

ACKNOWLEDGMENTS

This research was partially funded by the Australian Government through the Australian Research Council (project numbers DP210103092 and DP220103160). We acknowledge Grant PID2021-123611OB-I00 funded by MCIN/AEI/10.13039/501100011033 and by ERDF ''A way of making Europe''. P.A. Martorell acknowledges the support received from Universitat Politècnica de Catalunya and Santander Bank through an FPI fellowship (FPI-UPC 2019). This work was also supported by computational resources provided by the Australian Government through NCI under the National Computational Merit Allocation Scheme.